\documentclass[reqno, 12pt]{amsart}
 \usepackage[a4paper, hmargin = 1.25 in, vmargin = 1.5 in]{geometry}
 \usepackage{amsmath, amssymb, amsfonts}
 \usepackage{amsthm}
 
 \newtheorem{thm}{Theorem}[section]

 \newtheorem{lem}[thm]{Lemma}
 \newtheorem{prop}[thm]{Proposition}
 \theoremstyle{definition}
 \newtheorem{defn}{Definition}[section]
 \theoremstyle{remark}
 
 \numberwithin{equation}{section}

\begin{document}

\begin{center}
\begin{title}
\title{\bf\Large{{A Short Note on Lyapunov-type Inequalities for Hilfer Fractional Boundary Value Problems}}}
\end{title}

\vskip 0.25 in

\begin{author}
\author {Jagan Mohan Jonnalagadda\footnote[1]{Department of Mathematics, Birla Institute of Technology and Science Pilani, Hyderabad - 500078, Telangana, India. email: {j.jaganmohan@hotmail.com}}} 
\end{author}
\end{center}

\vskip 0.25 in

\noindent{\bf Abstract:} This paper deals with fractional boundary value problems involving the Hilfer fractional differential operator of order $1 < \alpha \leq 2$ and type $0 \leq \beta \leq 1$. We derive the corresponding Lyapunov-type inequalities for two prominent classes of Hilfer fractional boundary value problems (HFBVPs) involving separated and anti-periodic boundary conditions. For this purpose, we construct the associated Green's functions and deduce their important properties.

\vskip 0.25 in


\vskip 0.25 in


\vskip 0.25 in

\section{Introduction}
Fractional differential equations (FDEs) have proved to be valuable tools in modelling many phenomena in various fields of science and engineering. There has been a significant development in the theory and applications of FDEs in the last few decades; see the monographs of Podlubny \cite{P1}, Kilbas et al. \cite{KIL}, Hilfer \cite{H} and the references therein.

Lyapunov \cite{L} proved an inequality for Hill's equation associated with conjugate boundary conditions, known as Lyapunov's inequality. Cheng \cite{Ch} developed the discrete analogue of the inequality of Lyapunov for the first time. Many improvements in Lyapunov's inequality have been carried out due to its wide applications in oscillation theory, disconjugacy, eigenvalue problems etc. The study of Lyapunov-type inequalities for fractional boundary value problems has recently begun. In \cite{F1}, Ferreira replaced the second-order derivative with the $\alpha^{\text{th}}$-order, $\alpha \in (1, 2]$, Riemann--Liouville fractional derivative and obtained a Lyapunov-type inequality for the corresponding conjugate boundary value problem. Following this work, many authors gave valuable contributions to this topic. For an excellent introduction to the evolution of Lyapunov-type inequalities for ordinary differential equations, ordinary difference equations and fractional differential equations, we refer to \cite{Ag, N, N1, Pi, Ti} and the references therein. 

Recently, Ntouyas et al. \cite{N, N1} surveyed several generalizations of Lyapunov's inequality for fractional boundary value problems involving a variety of fractional derivative operators and boundary conditions. Motivated by these works, in this paper, we establish Lyapunov-type inequalities for the following fractional boundary value problems involving the Hilfer fractional differential operator of order $1 < \alpha \leq 2$ and type $0 \leq \beta \leq 1$. To our knowledge, no such work has yet been reported.
\begin{equation} \label{BVP A}
\begin{cases}
\left(D^{\alpha, \beta}_{a} y\right)(t) + q(t)y(t) = 0, \quad a < t < b,\\l\left(I_{a}^{(2-\alpha)(1-\beta)}y\right)(a)-m\left(DI_{a}^{(2-\alpha)(1-\beta)}y\right)(a)=0,\\n\left(I_{a}^{(2-\alpha)(1-\beta)}y\right)(b)+p\left(DI_{a}^{(2-\alpha)(1-\beta)}y\right)(b)=0,
\end{cases}
\end{equation}
and
\begin{equation} \label{BVP B}
\begin{cases}
\left(D^{\alpha, \beta}_{a} y\right)(t) + q(t)y(t) = 0, \quad a < t < b,\\\left(I_{a}^{(2-\alpha)(1-\beta)}y\right)(a)+\left(I_{a}^{(2-\alpha)(1-\beta)}y\right)(b)=0,\\\left(DI_{a}^{(2-\alpha)(1-\beta)}y\right)(a)+\left(DI_{a}^{(2-\alpha)(1-\beta)}y\right)(b)=0,
\end{cases}
\end{equation}
where $l$, $m$, $n$, $p$ are constants such that $l^2 + m^2 > 0$ and $n^2 + p^2 > 0$, $q : [a, b] \rightarrow \mathbb{R}$ is a continuous function, $D^{\alpha, \beta}_{a}$ denotes the Hilfer fractional differential operator of order $1 < \alpha \leq 2$ and type $0 \leq \beta \leq 1$ and $D$ denotes the first order differential operator.  

\section{Preliminaries}
In this section, we present some important definitions and results of fractional calculus which will be useful in the next section.

\begin{defn} \cite{KIL}
The Euler gamma function is given by $$\Gamma(z) = \int_0^\infty t^{z - 1} e^{-t} dt, \quad \Re(z) > 0.$$ Using the reduction formula $$ \Gamma (z + 1) = z\Gamma (z), \quad \Re(z) > 0,$$ the Euler gamma function can be extended to the half-plane $\Re(z) \leq 0$ except for $z \neq 0, -1, -2, \ldots$
\end{defn}

\begin{defn} \cite{KIL}
Let $a$, $b \in \mathbb{R}$ and $\alpha > 0$. The $\alpha^{\text{th}}$-order Riemann--Liouville fractional integral of a function $y : [a, b] \rightarrow \mathbb{R}$ is defined by 
\begin{equation}
\nonumber \left(I^{\alpha}_{a} y\right)(t) = \frac{1}{\Gamma(\alpha)} \int^{t}_{a} (t - s)^{\alpha - 1} y(s)ds, \quad t > a,
\end{equation}
provided the right-hand side exists. For $\alpha = 0$, we define $I^{\alpha}_{a}$ to be the identity map.
\end{defn}

\begin{defn} \cite{H1}
Let $a$, $b \in \mathbb{R}$, $\alpha > 0$, $0 \leq \beta \leq 1$ and choose $n \in \mathbb{N}$ such that $n - 1 < \alpha \leq n$. The $\alpha^{\text{th}}$-order and $\beta^{\text{th}}$-type  Hilfer fractional derivative of a function $y : [a, b] \rightarrow \mathbb{R}$ is defined by 
\begin{equation}
\nonumber (D^{\alpha, \beta}_{a}y)(t)=\left(I_{a}^{\beta(n-\alpha)}D^{n}I_{a}^{(n-\alpha)(1-\beta)}y\right)(t), \quad t > a,
\end{equation}
if the right-hand side exists. Here $D^{n} = \frac{d^{n}}{dt^{n}}$ denotes the classical $n^{\text{th}}$-order differential operator.
\end{defn}

\begin{defn} \cite{KIL}
We denote by $L(a, b)$ the space of Lebesgue measurable functions $y : [a, b] \rightarrow \mathbb{R}$ for which $$\|y\|_{L} = \int^{b}_{a}|y(t)|dt < \infty.$$
\end{defn}

\begin{defn} \cite{KIL}
We denote by $C[a, b]$ the space of continuous functions $y : [a, b] \rightarrow \mathbb{R}$ with the norm $$\|y\|_{C} = \max_{t \in [a, b]}|y(t)|.$$
\end{defn}

\begin{defn} \cite{KIL}
Let $AC[a, b]$ be the space of real valued functions $y$ which are absolutely continuous on $[a, b]$. We denote by $AC^n[a, b]$ the space of real valued functions $y$ which have continuous derivatives up to order $n - 1$ on $[a, b]$ such that $y^{(n - 1)} \in AC[a, b]$. In particular, $AC^1[a, b] = AC[a, b]$.
\end{defn}

\begin{defn} \cite{KIL}
Let $a$, $b \in \mathbb{R}$, $1 < \alpha \leq 2$, $0 \leq \beta \leq 1$ and $\gamma = (2 - \alpha)(1 - \beta) \in [0, 1)$. We denote by $C_{\gamma}[a, b]$ the weighted space of continuous functions $$C_{\gamma}[a, b] = \left\{y : (a, b] \rightarrow \mathbb{R} : (t - a)^{\gamma}y(t) \in C[a, b]\right\},$$ with the norm $$\|y\|_{C_{\gamma}} = \max_{t \in [a, b]}|(t - a)^{\gamma}y(t)|.$$
\end{defn}

\begin{lem} \cite{H1} \label{Power Rule}
For $a \in \mathbb{R}$, $\alpha > 0$, $0 \leq \beta \leq 1$ and $\mu > -1$, we have $$D^{\alpha, \beta}_{a} (t - a)^{\mu} = \frac{\Gamma(\mu+1)}{\Gamma(\mu+1-\alpha)}(t-a)^{\mu-\alpha}, \quad t > a.$$
\end{lem}

\begin{lem} \label{Solution 1}
For $a \in \mathbb{R}$, $1 < \alpha \leq 2$ and $0 \leq \beta \leq 1$, we have 
$$I^{(2-\alpha)(1-\beta)}_{a} (t - a)^{-(2-\alpha)(1-\beta)} = \Gamma(1-(2-\alpha)(1-\beta)), \quad t > a,$$ $$I^{(2-\alpha)(1-\beta)}_{a} (t - a)^{1-(2-\alpha)(1-\beta)} =(t-a)\Gamma(2-(2-\alpha)(1-\beta)), \quad t > a.$$ 
\end{lem}

\begin{lem} \cite{ZT} \label{Solution}
Let $y \in L(a, b)$, $n -1 < \alpha \leq n$, $n \in \mathbb{N}$, $\beta \in [0, 1]$ and $I_{a}^{(n - \alpha)(1 - \beta)}y \in AC^{k}[a, b]$. Then, we have
$$\left(I^{\alpha}_{a} D^{\alpha, \beta}_{a} y\right)(t) = y(t) - \sum^{n - 1}_{k = 0}\frac{(t-a)^{k-(n-\alpha)(1-\beta)}}{\Gamma(k-(n-\alpha)(1-\beta)+1)}\lim_{t\rightarrow a^{+}}\frac{d^k}{dt^{k}}\left(I_{a}^{(n-\alpha)(1-\beta)}y\right)(t),$$ for $t > a$.
\end{lem}

\begin{prop} \label{Max of f and g}
Let $f$ and $g$ be two nonnegative real valued functions defined on a set $S$. Assume that $f$ and $g$ attain their maximum values in $S$. Then, for each fixed $t \in S$, $$\vert{f(t) - g(t)} \vert \leq \max \left\{f(t), g(t)\right\} \leq \max \left\{\max_{t \in S}f(t), \max_{t \in S}g(t)\right\}.$$
\end{prop}

\section{Main Results}
In this section, we establish Lyapunov-type inequalities for the HFBVPs \eqref{BVP A} and \eqref{BVP B}. For this purpose, we consider the corresponding linear HFBVPs:
\begin{equation} \label{HIL FDE}
\begin{cases}
\left(D^{\alpha, \beta}_{a} y\right)(t) + h(t) = 0, \quad a < t < b,\\l\left(I_{a}^{(2-\alpha)(1-\beta)}y\right)(a)-m\left(DI_{a}^{(2-\alpha)(1-\beta)}y\right)(a)=0,\\n\left(I_{a}^{(2-\alpha)(1-\beta)}y\right)(b)+p\left(DI_{a}^{(2-\alpha)(1-\beta)}y\right)(b)=0,
\end{cases}
\end{equation}
and 
\begin{equation} \label{HIL-ANTI FDE}
\begin{cases}
\left(D^{\alpha, \beta}_{a} y\right)(t) + h(t) = 0, \quad a < t < b,\\ \left(I_{a}^{(2-\alpha)(1-\beta)}y\right)(a)+\left(I_{a}^{(2-\alpha)(1-\beta)}y\right)(b)=0,\\ \left(DI_{a}^{(2-\alpha)(1-\beta)}y\right)(a)+\left(DI_{a}^{(2-\alpha)(1-\beta)}y\right)(b)=0,
\end{cases}
\end{equation}
where $h : [a, b] \rightarrow \mathbb{R}$.

\begin{thm} \label{HIL Theorem}
Assume that $mn + lp + ln(b - a) \neq 0$. The unique solution of the HFBVP \eqref{HIL FDE} is given by
\begin{equation} \label{HIL Solution}
y(t) = \int^{b}_{a}G(t, s)h(s)ds, \quad a < t < b, 
\end{equation}
where
\begin{equation} \label{HIL Green}
G(t, s) = \begin{cases}
G_{1}(t,s), \quad a < s \leq t < b,\\
G_{2}(t,s), \quad a < t \leq s < b.
\end{cases}
\end{equation}
Here
$$G_{1}(t,s) = G_{2}(t, s) - \frac{(t - s)^{\alpha - 1}}{\Gamma(\alpha)},$$
and 
\begin{multline}
G_{2}(t,s) = (t - a)^{-(2 - \alpha)(1 - \beta)}(b - s)^{- \beta(2 - \alpha)} \\ \times \frac{\left[l(t - a) + m (\alpha - 1 + \beta (2 - \alpha))\right] \left[n(b - s) + p(1 - 2\beta + \alpha \beta)\right]}{\left[mn + lp + ln(b - a)\right]\Gamma(2 - 2\beta + \alpha \beta)\Gamma(2 - (2 - \alpha)(1 - \beta))}.
\end{multline}
\end{thm}

\begin{proof}
Applying the $\alpha^{\text{th}}$-order Riemann--Liouville fractional integration operator on both sides of \eqref{HIL FDE} and applying the Lemma \ref{Solution}, we get
\begin{equation} \label{Sol 1}
y(t) = C_{1}\frac{(t-a)^{-(2-\alpha)(1-\beta)}}{\Gamma(1-(2-\alpha)(1-\beta))} + C_{2}\frac{(t-a)^{1-(2-\alpha)(1-\beta)}}{\Gamma(2-(2-\alpha)(1-\beta))}-\int^{t}_{a}\frac{(t-s)^{\alpha-1}}{\Gamma(\alpha)}h(s)ds.
\end{equation}
Now, using Lemma \ref{Solution 1} to \eqref{Sol 1}, we have \begin{equation} \label{Sol 2}
\left(I_{a}^{(2-\alpha)(1-\beta)}y\right)(t)= C_1+C_2(t-a)-\int^{t}_{a}\frac{(t-s)^{1-2\beta+\alpha\beta}}{\Gamma(2-2\beta+\alpha\beta)}h(s)ds
\end{equation}
and
\begin{equation} \label{Sol 3}
\left(DI_{a}^{(2-\alpha)(1-\beta)}y\right)(t)= C_2-\int^{t}_{a}\frac{(t-s)^{-2\beta+\alpha\beta}}{\Gamma(1-2\beta+\alpha\beta)}h(s)ds.
\end{equation}
Using the boundary conditions $$l\left(I_{a}^{(2-\alpha)(1-\beta)}y\right)(a)-m\left(DI_{a}^{(2-\alpha)(1-\beta)}y\right)(a)=0$$ and $$n\left(I_{a}^{(2-\alpha)(1-\beta)}y\right)(b)+p\left(DI_{a}^{(2-\alpha)(1-\beta)}y\right)(b)=0$$ to \eqref{Sol 2} and \eqref{Sol 3}, we get 
\begin{multline}
\nonumber C_{1} = \frac{m}{\left[mn + lp + ln(b - a)\right]\Gamma(2 - 2 \beta + \alpha \beta)} \\ \times \int^{b}_{a} (b - s)^{- \beta(2 - \alpha)}\left[n(b - s) + p(1 - 2\beta + \alpha \beta)\right]h(s)ds.
\end{multline}
and 
\begin{multline}
\nonumber C_{2} = \frac{l}{\left[mn + lp + ln(b - a)\right]\Gamma(2 - 2 \beta + \alpha \beta)} \\ \times \int^{b}_{a} (b - s)^{- \beta(2 - \alpha)}\left[n(b - s) + p(1 - 2\beta + \alpha \beta)\right]h(s)ds.
\end{multline}
Now, substituting the values of $C_1$ and $C_2$ in \eqref{Sol 1} and rearranging the terms, we obtain \eqref{HIL Solution}. The proof is complete.
\end{proof}

\begin{thm} \label{HIL-ANTI Theorem}
The unique solution of the HFBVP \eqref{HIL-ANTI FDE} is given by
\begin{equation} \label{HIL-ANTI Solution}
y(t) = \int^{b}_{a}\bar{G}(t, s)h(s)ds, \quad a < t < b,
\end{equation}
where
\begin{equation} \label{HIL-ANTI Green}
\bar{G}(t, s) = \begin{cases}
\bar{G}_{1}(t,s), \quad a < s \leq t < b,\\
\bar{G}_{2}(t,s), \quad a < t \leq s < b.
\end{cases}
\end{equation}
Here $$\bar{G}_{1}(t,s)=\bar{G}_{2}(t,s)- \frac{(t-s)^{\alpha-1}}{\Gamma(\alpha)},$$
and $$\bar{G}_{2}(t,s)=\frac{(t-a)^{-(2-\alpha)(1-\beta)}(b-s)^{-\beta(2-\alpha)}\left[\frac{(t-a)}{(\alpha + \beta - \alpha \beta)}-\frac{(b-s)}{(1 - \beta(2 - \alpha))}-\frac{(b-a)}{2}\right]}{2 \Gamma(1-2\beta+\alpha \beta)\Gamma(1-(2-\alpha)(1-\beta))}.$$
\end{thm}
 
\begin{proof}
Using the boundary conditions $$\left(I_{a}^{(2-\alpha)(1-\beta)}y\right)(a)+\left(I_{a}^{(2-\alpha)(1-\beta)}y\right)(b)=0$$ and $$\left(DI_{a}^{(2-\alpha)(1-\beta)}y\right)(a)+\left(DI_{a}^{(2-\alpha)(1-\beta)}y\right)(b)=0$$ to \eqref{Sol 2} and \eqref{Sol 3}, we get $$C_1 = - \frac{(b - a)}{4}\int^{b}_{a}\frac{(b-s)^{-2\beta+\alpha\beta}}{\Gamma(1-2\beta+\alpha\beta)}h(s)ds - \frac{1}{2} \int^{b}_{a}\frac{(b-s)^{1-2\beta+\alpha\beta}}{\Gamma(2-2\beta+\alpha\beta)}h(s)ds$$ and $$C_{2} = \frac{1}{2} \int^{b}_{a}\frac{(b-s)^{-2\beta+\alpha\beta}}{\Gamma(1-2\beta+\alpha\beta)}h(s)ds.$$
Now, substituting the values of $C_1$ and $C_2$ in \eqref{Sol 1} and rearranging the terms, we obtain \eqref{HIL-ANTI Solution}. The proof is complete.
\end{proof}

\begin{thm} \label{HIL Green Max}
Assume that $l$, $m$, $n$, $p \geq 0$ and $mn + lp + ln(b - a) > 0$. Denote by 
\begin{align*}
H(t, s) & = (t - a)^{(2 - \alpha)(1 - \beta)}(b - s)^{\beta(2 - \alpha)}G(t, s) \\ & = \begin{cases}
(t - a)^{(2 - \alpha)(1 - \beta)}(b - s)^{\beta(2 - \alpha)}G_{1}(t,s), \quad a \leq s \leq t \leq b,\\
(t - a)^{(2 - \alpha)(1 - \beta)}(b - s)^{\beta(2 - \alpha)}G_{2}(t,s), \quad a \leq t \leq s \leq b. 
\end{cases} \\
& = \begin{cases}
H_{1}(t,s), \quad a \leq s \leq t \leq b,\\
H_{2}(t,s), \quad a \leq t \leq s \leq b.
\end{cases}
\end{align*}
Then, $$\left|H(t, s)\right| \leq \max \left\{\Omega, \frac{(b - a)}{\Gamma(\alpha)}\right\}, \quad (t, s) \in [a, b] \times [a, b],$$ where $$\Omega = \frac{\left[l(b - a) + m (\alpha - 1 + \beta (2 - \alpha))\right] \left[n(b - a) + p(1 - 2\beta + \alpha \beta)\right]}{\left[mn + lp + ln(b - a)\right]\Gamma(2 - 2\beta + \alpha \beta)\Gamma(2 - (2 - \alpha)(1 - \beta))}.$$
\end{thm}

\begin{proof}
We have $$H_1(t,s) = H_2(t, s) - \frac{(t - s)^{\alpha - 1}(t - a)^{(2 - \alpha)(1 - \beta)}(b - s)^{\beta(2 - \alpha)}}{\Gamma(\alpha)}, \quad a \leq s \leq t \leq b,$$ and $$H_2(t, s) = \frac{\left[l(t - a) + m (\alpha - 1 + \beta (2 - \alpha))\right] \left[n(b - s) + p(1 - 2\beta + \alpha \beta)\right]}{\left[mn + lp + ln(b - a)\right]\Gamma(2 - 2\beta + \alpha \beta)\Gamma(2 - (2 - \alpha)(1 - \beta))},$$ for $a \leq t \leq s \leq b$. Clearly, $$H_2(t, s) \geq 0, \quad (t, s) \in [a, b] \times [a, b],$$ and $$\frac{(t - s)^{\alpha - 1}(t - a)^{(2 - \alpha)(1 - \beta)}(b - s)^{\beta(2 - \alpha)}}{\Gamma(\alpha)} \geq 0, \quad a \leq s \leq t \leq b.$$ Now, we apply Proposition \ref{Max of f and g}. For $a \leq s \leq t \leq b$, we obtain
\begin{align*}
& \left|H_1(t, s)\right| = \left|H_2(t,s) - \frac{(t - s)^{\alpha - 1}(t - a)^{(2 - \alpha)(1 - \beta)}(b - s)^{\beta(2 - \alpha)}}{\Gamma(\alpha)}\right| \\ \leq & \max \left\{\max_{a \leq s \leq t \leq b}H_2(t, s), \max_{a \leq s \leq t \leq b}\left[\frac{(t - s)^{\alpha - 1}(t - a)^{(2 - \alpha)(1 - \beta)}(b - s)^{\beta(2 - \alpha)}}{\Gamma(\alpha)}\right]\right\}.
\end{align*}
Denote by $$K(t, s) = \frac{(t - s)^{\alpha - 1}(t - a)^{(2 - \alpha)(1 - \beta)}(b - s)^{\beta(2 - \alpha)}}{\Gamma(\alpha)}, \quad a \leq s \leq t \leq b.$$ For a fixed $t \in [a, b]$, consider $$\frac{\partial}{\partial s}K(t, s) = \frac{\beta(2 - \alpha)(t - s)^{\alpha - 1}(t - a)^{(2 - \alpha)(1 - \beta)}(b - s)^{\beta(2 - \alpha) - 1}}{\Gamma(\alpha)}.$$ Since $1 < \alpha \leq 2$, $0 \leq \beta \leq 1$, $a \leq s \leq t \leq b$, we have $$\frac{\partial}{\partial s}K(t, s) > 0, \quad \quad a \leq s \leq t \leq b,$$ implying that $K(t, s)$ is a decreasing function of $s$ for a fixed $t \in [a, b]$. Then, $$K(t, s) \leq K(t, a), \quad a \leq s \leq t \leq b.$$ Denote by $$f(t) = K(t, a) = \frac{(t - a)^{\alpha - 1}(t - a)^{(2 - \alpha)(1 - \beta)}(b - a)^{\beta(2 - \alpha)}}{\Gamma(\alpha)}, \quad a \leq t \leq b.$$ Consider $$f'(t) = \frac{(1 - \beta (2 - \alpha))(t - a)^{- \beta (2 - \alpha)}(b - a)^{\beta(2 - \alpha)}}{\Gamma(\alpha)}, \quad a \leq t \leq b.$$ Since $1 < \alpha \leq 2$, $0 \leq \beta \leq 1$, $a \leq t \leq b$, we have $$f'(t) > 0, \quad a \leq t \leq b,$$ implying that $f$ is an increasing function of $t$ for $t \in [a, b]$. Thus, $$f(t) \leq f(b), \quad a \leq t \leq b.$$ Therefore, $$K(t, s) \leq K(b, a), \quad a \leq s \leq t \leq b.$$ That is, $$\max_{a \leq s \leq t \leq b} K(t, s) = \frac{(b - a)}{\Gamma(\alpha)}.$$ We also know that $$\max_{(t, s) \in [a, b] \times [a, b]} H_2(t, s) = \Omega.$$ Therefore, $$\left|H_1(t, s)\right| \leq \max \left\{\Omega, \frac{(b - a)}{\Gamma(\alpha)}\right\}, \quad a \leq s \leq t \leq b.$$ Hence, for all $(t, s) \in [a, b] \times [a, b]$, 
\begin{align*}
& \left|H(t, s)\right| = \max \left\{\max_{a \leq t \leq s \leq b}\left|H_2(t, s)\right|, \max_{a \leq s \leq t \leq b}\left|H_1(t, s)\right|\right\} \\ & = \max \left\{\Omega, \frac{(b - a)}{\Gamma(\alpha)}\right\}.
\end{align*}
The proof is complete.
\end{proof}

\begin{thm} \label{HIL Green Max 1}
Denote by 
\begin{align*}
\bar{H}(t, s) & = (t - a)^{(2 - \alpha)(1 - \beta)}(b - s)^{\beta(2 - \alpha)}\bar{G}(t, s) \\ & = \begin{cases}
(t - a)^{(2 - \alpha)(1 - \beta)}(b - s)^{\beta(2 - \alpha)}\bar{G}_{1}(t,s), \quad a \leq s \leq t \leq b,\\
(t - a)^{(2 - \alpha)(1 - \beta)}(b - s)^{\beta(2 - \alpha)}\bar{G}_{2}(t,s), \quad a \leq t \leq s \leq b. 
\end{cases} \\
& = \begin{cases}
\bar{H}_{1}(t,s), \quad a \leq s \leq t \leq b,\\
\bar{H}_{2}(t,s), \quad a \leq t \leq s \leq b.
\end{cases}
\end{align*}
Then, $$\left|\bar{H}(t, s)\right| \leq \frac{(b - a)}{A} \left[\frac{1}{(1 - \beta(2 - \alpha))} + \frac{1}{2}\right] + \frac{(b - a)}{\Gamma(\alpha)}, \quad (t, s) \in [a, b] \times [a, b],$$ where $$A = 2 \Gamma(1 - 2\beta + \alpha \beta)\Gamma(1 - (2 - \alpha)(1 - \beta)).$$ 
\end{thm}

\begin{proof}
We have $$\bar{H}_1(t,s) = \bar{H}_2(t, s) - \frac{(t - s)^{\alpha - 1}(t - a)^{(2 - \alpha)(1 - \beta)}(b - s)^{\beta(2 - \alpha)}}{\Gamma(\alpha)}, \quad a \leq s \leq t \leq b,$$ and $$\bar{H}_2(t, s) = \frac{1}{A}\left[\frac{(t - a)}{(\alpha + \beta - \alpha \beta)} - \frac{(b - s)}{(1 - \beta(2 - \alpha))} - \frac{(b - a)}{2}\right], \quad a \leq t \leq s \leq b.$$ Clearly, $A > 0$ and $$\frac{(t - s)^{\alpha - 1}(t - a)^{(2 - \alpha)(1 - \beta)}(b - s)^{\beta(2 - \alpha)}}{\Gamma(\alpha)} \geq 0, \quad a \leq s \leq t \leq b.$$ For a fixed $s \in [a, b]$, we have 
\begin{equation} \label{Max - 1}
\frac{\partial}{\partial t} \bar{H}_2(t, s) = \frac{1}{A(\alpha + \beta - \alpha \beta)} > 0,
\end{equation}
implying that $\bar{H}_2(t, s)$ is an increasing function of $t$. Thus, we have $$\max_{a \leq t \leq s}\left|\bar{H}_2(t, s)\right| = \max \left\{\left|\bar{H}_2(a, s)\right|, \left|\bar{H}_2(s, s)\right|\right\}.$$ We observe that $\bar{H}_2(s, s)$ is an increasing function of $s$, since $$\frac{d}{ds} \bar{H}_2(s, s) = \frac{1}{A}\left[\frac{1}{(\alpha + \beta - \alpha \beta)} + \frac{1}{(1 - \beta(2 - \alpha))}\right] > 0.$$ Therefore, we have 
\begin{align*}
& \max_{a \leq s \leq b}\bar{H}_2(s, s) \\ & = \max \left\{\left|\bar{H}_2(a, a)\right|, \left|\bar{H}_2(b, b)\right|\right\} \\ & =  \frac{1}{A} \max \left\{\left|- \frac{(b - a)}{(1 - \beta(2 - \alpha))} - \frac{(b - a)}{2}\right|, \left|\frac{(b - a)}{(\alpha + \beta - \alpha \beta)} - \frac{(b - a)}{2}\right|\right\} \\ & =  \frac{(b - a)}{A} \max \left\{\left[\frac{1}{(1 - \beta(2 - \alpha))} + \frac{1}{2}\right], \left|\frac{1}{(\alpha + \beta - \alpha \beta)} - \frac{1}{2}\right|\right\}.
\end{align*}
Since $1 < \alpha \leq 2$ and $0 \leq \beta \leq 1$, we have $\alpha + \beta - \alpha \beta < 2$ implying that 
\begin{equation} \label{Max 0}
\frac{1}{(\alpha + \beta - \alpha \beta)} > \frac{1}{2}.
\end{equation}
So, $$\max_{a \leq s \leq b}\bar{H}_2(s, s) = \frac{(b - a)}{A} \max \left\{\left[\frac{1}{(1 - \beta(2 - \alpha))} + \frac{1}{2}\right], \left[\frac{1}{(\alpha + \beta - \alpha \beta)} - \frac{1}{2}\right]\right\}.$$
Now, consider 
\begin{align*}
& \frac{1}{(1 - \beta(2 - \alpha))} - \frac{1}{(\alpha + \beta - \alpha \beta)} \\ & = \frac{1}{(1 - \beta(2 - \alpha))} - \frac{1}{(\alpha + \beta - \alpha \beta)} \\ & = \frac{(\alpha + \beta - \alpha \beta) - (1 - \beta(2 - \alpha))}{(1 - \beta(2 - \alpha))(\alpha + \beta - \alpha \beta)} \\ & =  \frac{(\alpha + 3 \beta - 2 \alpha \beta - 1)}{(1 - \beta(2 - \alpha))(\alpha + \beta - \alpha \beta)}.
\end{align*}
Since $1 < \alpha \leq 2$ and $0 \leq \beta \leq 1$, we have $1 - \beta(2 - \alpha) > 0$, $\alpha + \beta - \alpha \beta > 0$ and $$\alpha + 3 \beta - 2 \alpha \beta - 1 = (\alpha - 1)(1 - \beta) + \beta (2 - \alpha) \geq 0,$$ implying that 
\begin{equation} \label{Max 1}
\frac{1}{(1 - \beta(2 - \alpha))} > \frac{1}{(\alpha + \beta - \alpha \beta)}.
\end{equation}
That is $$\left[\frac{1}{(1 - \beta(2 - \alpha))} + \frac{1}{2}\right] > \left[\frac{1}{(\alpha + \beta - \alpha \beta)} - \frac{1}{2}\right].$$ Therefore, 
\begin{equation} \label{Max 2}
\max_{a \leq s \leq b}\bar{H}_2(s, s) = \frac{(b - a)}{A}\left[\frac{1}{(1 - \beta(2 - \alpha))} + \frac{1}{2}\right].
\end{equation}
Now, we consider $$\frac{d}{ds} \bar{H}_2(a, s) = \frac{1}{A}\left[\frac{1}{(1 - \beta(2 - \alpha))}\right] > 0.$$ So, $\bar{H}_2(a, s)$ is an increasing function of $s$. Thus, we have
\begin{align*}
\max_{a \leq s \leq b}\bar{H}_2(a, s) & = \max \left\{\left|\bar{H}_2(a, a)\right|, \left|\bar{H}_2(a, b)\right|\right\} \\ & =  \frac{1}{A} \max \left\{\left|- \frac{(b - a)}{(1 - \beta(2 - \alpha))} - \frac{(b - a)}{2}\right|, \left|- \frac{(b - a)}{2}\right|\right\} \\ & = \frac{(b - a)}{A} \max \left\{\left[\frac{1}{(1 - \beta(2 - \alpha))} + \frac{1}{2}\right], \frac{1}{2}\right\} \\ & = \frac{(b - a)}{A} \left[\frac{1}{(1 - \beta(2 - \alpha))} + \frac{1}{2}\right].
\end{align*}
Hence, we have $$\max_{a \leq t \leq s \leq b}\left|\bar{H}_2(t, s)\right| = \frac{(b - a)}{A} \left[\frac{1}{(1 - \beta(2 - \alpha))} + \frac{1}{2}\right].$$ Now, for $a \leq s \leq t \leq b$, consider $\bar{H}_2(t, s)$. For a fixed $s \in [a, b]$, it follows from \eqref{Max - 1} that $\bar{H}_2(t, s)$ is an increasing function of $t$. Thus, we have $$\max_{s \leq t \leq b}\left|\bar{H}_2(t, s)\right| = \max \left\{\left|\bar{H}_2(b, s)\right|, \left|\bar{H}_2(s, s)\right|\right\}.$$ Consider $$\frac{d}{ds} \bar{H}_2(b, s) = \frac{1}{A}\left[\frac{1}{(1 - \beta(2 - \alpha))}\right] > 0.$$ So, $\bar{H}_2(b, s)$ is an increasing function of $s$. Thus, we have
\begin{align*}
\max_{a \leq s \leq b}\bar{H}_2(b, s) & = \max \left\{\left|\bar{H}_2(b, a)\right|, \left|\bar{H}_2(b, b)\right|\right\} \\ & = \frac{1}{A} \max \Bigg{\{}\left|\frac{(b - a)}{(\alpha + \beta - \alpha \beta)} - \frac{(b - a)}{(1 - \beta(2 - \alpha))} - \frac{(b - a)}{2}\right|, \\ & \quad \left|\frac{(b - a)}{(\alpha + \beta - \alpha \beta)} - \frac{(b - a)}{2}\right|\Bigg{\}} \\ & = \frac{(b - a)}{A} \max \Bigg{\{}\left|\frac{1}{(\alpha + \beta - \alpha \beta)} - \frac{1}{(1 - \beta(2 - \alpha))} - \frac{1}{2}\right|, \\ & \quad \left|\frac{1}{(\alpha + \beta - \alpha \beta)} - \frac{1}{2}\right|\Bigg{\}} \\ & = \frac{(b - a)}{A} \max \Bigg{\{}\left[\frac{1}{2} + \frac{1}{(1 - \beta(2 - \alpha))} - \frac{1}{(\alpha + \beta - \alpha \beta)}\right], \\ & \quad \left[\frac{1}{(\alpha + \beta - \alpha \beta)} - \frac{1}{2}\right]\Bigg{\}}. \quad (\text{By \eqref{Max 0}, \eqref{Max 1}})
\end{align*}
Then, it follows from \eqref{Max 2} that $$\max_{a \leq s \leq t \leq b}\left|\bar{H}_2(t, s)\right| = \frac{(b - a)}{A} \left[\frac{1}{(1 - \beta(2 - \alpha))} + \frac{1}{2}\right].$$ Denote by $$K(t, s) = \frac{(t - s)^{\alpha - 1}(t - a)^{(2 - \alpha)(1 - \beta)}(b - s)^{\beta(2 - \alpha)}}{\Gamma(\alpha)}, \quad a \leq s \leq t \leq b.$$ It follows from the proof of the above theorem, we obtain $$\max_{a \leq s \leq t \leq b} K(t, s) = \frac{(b - a)}{\Gamma(\alpha)}.$$ Now,  for $a \leq s \leq t \leq b$, consider
\begin{align*}
\left|\bar{H}_1(t, s)\right| & = \left|\bar{H}_2(t, s) - K(t, s)\right| \\ & \leq \left|\bar{H}_2(t, s)\right| + \left|K(t, s)\right| \\ & \leq \frac{(b - a)}{A} \left[\frac{1}{(1 - \beta(2 - \alpha))} + \frac{1}{2}\right] + \frac{(b - a)}{\Gamma(\alpha)}.
\end{align*}
Therefore, for all $(t, s) \in [a, b] \times [a, b]$, 
\begin{align*}
& \left|\bar{H}(t, s)\right| = \max \left\{\max_{a \leq t \leq s \leq b}\left|\bar{H}_2(t, s)\right|, \max_{a \leq s \leq t \leq b}\left|\bar{H}_1(t, s)\right|\right\} \\ & = \frac{(b - a)}{A} \left[\frac{1}{(1 - \beta(2 - \alpha))} + \frac{1}{2}\right] + \frac{(b - a)}{\Gamma(\alpha)}.
\end{align*}
The proof is complete.
\end{proof}

Now, we are able to develop Lyapunov-type inequalities for the HFBVPs \eqref{BVP A} and \eqref{BVP B}.

\begin{thm} \label{HIL Lyp}
Assume that $l$, $m$, $n$, $p \geq 0$ and $mn + lp + ln(b - a) > 0$. If the HFBVP \eqref{BVP A} has a nontrivial solution, then
\begin{equation}
\int^{b}_{a}(s - a)^{(\alpha - 2)(1 - \beta)}(b - s)^{\beta(\alpha - 2)}\big{|}q(s)\big{|}ds \geq \frac{1}{\Lambda},
\end{equation}
where $$\Lambda = \max \left\{\Omega, \frac{(b - a)}{\Gamma(\alpha)}\right\}.$$
\end{thm}

\begin{proof}
We know that $\mathfrak{B} = C_{(2 - \alpha)(1 - \beta)}[a, b]$ is a Banach space with the norm $$\|y\|_{C_{(2 - \alpha)(1 - \beta)}} = \max_{t \in [a, b]}\Big{\vert}{(t - a)^{(2 - \alpha)(1 - \beta)}y(t)}\Big{\vert}.$$ It follows from Theorem \ref{HIL Theorem} that $y$ is a solution of \eqref{BVP A} if and only if $y$ is a solution of the Fredholm integral equation
\begin{equation*}
y(t) = \int^{b}_{a}G(t, s)q(s)y(s)ds.
\end{equation*}
Consider 
\begin{align*}
& \|y\|_{C_{(2 - \alpha)(1 - \beta)}} \\ & = \max_{t \in [a, b]}\Big{\vert}{(t - a)^{(2 - \alpha)(1 - \beta)}y(t)}\Big{\vert} \\ & = \max_{t \in [a, b]}\left|(t - a)^{(2 - \alpha)(1 - \beta)}\int^{b}_{a}G(t, s)q(s)y(s)ds\right| \\ & \leq \max_{t \in [a, b]}\Big{[}(t - a)^{(2 - \alpha)(1 - \beta)}\int^{b}_{a}\vert{G(t, s)}\vert\big{|}q(s)\big{|}\big{|}y(s)\big{|}ds\Big{]} \\ & = \max_{t \in [a, b]}\Big{[}(t - a)^{(2 - \alpha)(1 - \beta)}\int^{b}_{a}(s - a)^{- (2 - \alpha)(1 - \beta)} \vert{G(t, s)}\vert \big{|}q(s)\big{|} \\ & \quad \left[(s - a)^{(2 - \alpha)(1 - \beta)} \big{|}y(s)\big{|}\right]ds\Big{]} \\ & \leq \|y\|_{C_{(2 - \alpha)(1 - \beta)}} \max_{t \in [a, b]} \Big{[}(t - a)^{(2 - \alpha)(1 - \beta)} \int^{b}_{a}(s - a)^{- (2 - \alpha)(1 - \beta)} \vert{G(t, s)}\vert \big{|}q(s)\big{|}ds\Big{]} \\ & = \|y\|_{C_{(2 - \alpha)(1 - \beta)}} \int^{b}_{a}\max_{t \in [a, b]} \Big{|}(t - a)^{(2 - \alpha)(1 - \beta)}(b - s)^{\beta (2 - \alpha)}G(t, s)\Big{|} \\ & \quad (s - a)^{- (2 - \alpha)(1 - \beta)}(b - s)^{- \beta (2 - \alpha)} \big{|}q(s)\big{|}ds  \\ & = \|y\|_{C_{(2 - \alpha)(1 - \beta)}} \int^{b}_{a}\left[\max_{t \in [a, b]} \left|H(t, s)\right|\right](s - a)^{- (2 - \alpha)(1 - \beta)}(b - s)^{- \beta (2 - \alpha)} \big{|}q(s)\big{|}ds \\ & \leq \Lambda \|y\|_{C_{(2 - \alpha)(1 - \beta)}} \int^{b}_{a}(s - a)^{- (2 - \alpha)(1 - \beta)}(b - s)^{- \beta (2 - \alpha)} \big{|}q(s)\big{|}ds
\end{align*} implying that 
$$1 \leq \Lambda \int^{b}_{a}(s - a)^{- (2 - \alpha)(1 - \beta)}(b - s)^{- \beta (2 - \alpha)} \big{|}q(s)\big{|}ds.$$ The proof is complete.
\end{proof}

\begin{thm} \label{HIL Lyp 1}
If the HFBVP \eqref{BVP B} has a nontrivial solution, then
\begin{equation}
\int^{b}_{a}(s - a)^{(\alpha - 2)(1 - \beta)}(b - s)^{\beta(\alpha - 2)}\big{|}q(s)\big{|}ds \geq \frac{1}{\Theta},
\end{equation}
where $$\Theta = \frac{(b - a)}{A} \left[\frac{1}{(1 - \beta(2 - \alpha))} + \frac{1}{2}\right] + \frac{(b - a)}{\Gamma(\alpha)}.$$
\end{thm}

\begin{proof}
The proof is similar to the proof of Theorem \ref{HIL Lyp}. So, we omit it.
\end{proof}

\section*{Conclusion}
In this article, we derived the corresponding Lyapunov-type inequalities for two prominent classes of HFBVPs \eqref{BVP A} and \eqref{BVP B} involving separated and anti-periodic boundary conditions, respectively. For this purpose, we constructed the associated Green's functions and deduced their important properties.

\end{document}